\documentclass{journal}
\usepackage{fullpage}
\usepackage{generic}
\usepackage{cite}
\usepackage{amsmath,amssymb,amsfonts}
\usepackage{algorithmicx}
\usepackage{algorithm}
\usepackage{algpseudocode}
\usepackage{graphicx}
\usepackage{textcomp}
\usepackage{url}
\usepackage{float}
\usepackage{textcomp}
\usepackage{tabularx,booktabs,makecell}
\usepackage{setspace}
\usepackage{subcaption}
\usepackage{pifont}
\usepackage{nomencl}
\usepackage{etoolbox}
\def\BibTeX{{\rm B\kern-.05em{\sc i\kern-.025em b}\kern-.08em
T\kern-.1667em\lower.7ex\hbox{E}\kern-.125emX}}

\algnewcommand{\CallSolve}{\Call{SolveOptimalScheduling}}
\algnewcommand{\CallSolveMILP}{\Call{SolveMILP}}

\renewcommand\nomgroup[1]{%
  \item[\bfseries
  \ifstrequal{#1}{A}{\textit{Abbreviations}}{%
  \ifstrequal{#1}{B}{\textit{Sets and Indices}}{%
  \ifstrequal{#1}{C}{\textit{Parameters}}{%
  \ifstrequal{#1}{D}{\textit{Decision Variables}}{%
  \ifstrequal{#1}{E}{\textit{Derived or Intermediate Quantities}}{}}}}}]%
}

\renewcommand{\arraystretch}{1.3}
\newcolumntype{Y}{>{\raggedright\arraybackslash}X}
\newcommand{\cmark}{\ding{51}}  
\newcommand{\xmark}{\ding{55}}

\begin{document}
\title{Dynamic Menu-Based Pricing for Electric Vehicle Charging with Vehicle-to-Grid Integration}
\author{
Mozhdeh Hematiboroujeni,
Pierre Le Bodic,
Adel N. Toosi,
and \\ Markus Wagner
\thanks{Mozhdeh Hematiboroujeni, Pierre Le Bodic, and Markus Wagner are with the Department of Information Technology, Monash University, Melbourne, VIC 3800, Australia.\\
Adel N. Toosi is with the School of Computing and Information Systems, The University of Melbourne, Australia.\\
Emails: Mozhdeh.Hematiboroujeni@monash.edu, Pierre.Lebodic@monash.edu, adel.toosi@unimelb.edu.au, Markus.Wagner@monash.edu}
}
\maketitle
\begin{abstract}
The number of electric vehicles is rapidly increasing worldwide. This growth brings significant environmental benefits but also introduces new challenges: uncoordinated charging can place stress on the grid, particularly during peak hours. Beyond these challenges lies the opportunity for electric vehicles to feed energy back to the grid (V2G), which helps balance supply and demand and supports renewable energy. However, current pricing schemes such as time-of-use tariffs provide little incentive for discharging. To study incentive design in a realistic context, we focus on a parking lot operator who manages multiple EV chargers. We propose a menu-based pricing mechanism in which each EV declares its energy requirement and parking duration; given the retail real-time electricity prices, the operator offers a menu of options that trade off the allowed level of discharging and the associated price. We formulate this interaction as a bilevel optimization problem and reformulate it into a single-level model. 
Results show that, relative to a no-V2G baseline, the proposed mechanism increases operator profit by 30\% and reduces EV payments by 17\%. Compared to widely used tariff baselines, it improves operator profit by 22–29\%, lowers EV payments by 9–18\% and increases V2G contribution by 87–235\%.
Overall, the results show that the proposed dynamic menu-based pricing framework provides a practical, computationally efficient, and economically advantageous approach for real-time EV charging and V2G integration.

\end{abstract}

\renewcommand{\bottomfraction}{0.999}
\renewcommand{\topfraction}{0.999}
\renewcommand{\dbltopfraction}{0.999}
\setcounter{bottomnumber}{9}
\setcounter{topnumber}{9}
\setcounter{dbltopnumber}{9}
\setcounter{totalnumber}{9}
\section{Introduction}
\label{sec:intro}
The increasing adoption of electric vehicles (EVs) provides a promising opportunity for enhancing grid flexibility and supporting renewable energy integration~\cite{Barman2023RenewableEnergyIntegration}. This rapid growth, while promising for environmental goals, poses significant challenges to the power grid, particularly at parking lots equipped with EV chargers, where uncoordinated charging during peak demand periods can strain distribution networks~\cite{Muratori2019}. Therefore, effective coordination mechanisms are essential to manage EV charging while maintaining grid reliability and user convenience.

Vehicle-to-grid (V2G) technology offers a powerful solution to these challenges by enabling bidirectional energy exchange between EVs and the grid~\cite{Kempton2005}. Through coordinated charging and discharging, V2G-equipped EVs can provide valuable grid services such as peak load shaving, ancillary services~\cite{KarimiArpanahi2024BatteryScheduling}, and renewable integration~\cite{Noel2019}. However, the practical deployment of V2G has been limited by both technical and behavioral barriers. From the user’s perspective, concerns about battery degradation and inadequate financial rewards often discourage participation~\cite{Parsons2014,Sagaria2025V2GCompensation}. Moreover, current charging infrastructure and pricing mechanisms, such as time-of-use (ToU) tariffs, lack the flexibility to adapt to real-time grid conditions or incentivize discharging during periods of high demand~\cite{Sortomme2011, Hu2016}. These limitations highlight the need for advanced pricing strategies that align the interests of EV owners, parking lot operators, and grid while respecting distribution network constraints~\cite{ghosh2018v2g},~\cite{Nimalsiri2023}.

To address these challenges, previous research has explored a range of coordination and pricing strategies for EV charging, including dynamic pricing that adjusts electricity rates based on grid or market conditions~\cite{Silva2025CarbonAware}, demand response programs that encourage EV owners to modify their charging behavior in response to price signals or control requests~\cite{Kakkar2024EVDRSurvey}, advanced scheduling methods that jointly optimize charging decisions to reduce system costs and network stress~\cite{Meng2025FourStage}. While these approaches have provided valuable insights, they often rely on centralized optimization or assume perfect user information, making them less suitable for large-scale, user-driven V2G environments.
Among the emerging alternatives, menu-based pricing has shown strong potential to bridge the gap between economic efficiency and user flexibility~\cite{GhoshAggarwal2018Menu,Lu2023DeadlineMenu}.
In this mechanism, the operator offers each EV a menu of options for charging, allowing the EV to self-select the option that best matches its preferences. However, most existing menu-based studies simplify the problem by assuming static electricity tariffs or one-way (charging-only) energy flow. Moreover, they typically neglect real-time wholesale price variation, feeder capacity limits, and the sequential nature of EV arrivals, all of which are critical for real-world implementation.

In this context, our work proposes a novel dynamic menu-based pricing mechanism to enhance the management of EV charging with V2G capabilities in a parking lot setting. The proposed model introduces a real-time pricing strategy that provides each EV with a personalized menu of charging and discharging options tailored to its energy requirements and parking duration. By designing incentives that encourage V2G participation, our mechanism improves economic outcomes for parking lot operators and EV owners, while also providing support to the grid during peak demand periods. This framework is developed and validated through a comprehensive case study, demonstrating its applicability in a simulated urban setting.

Our contribution lies in combining real-time, per-EV menu updates with V2G support and feeder constraints, which are features that no prior menu-based pricing work has addressed jointly. Building on this foundation, the key contributions of this paper are as follows:

\begin{enumerate}
    \item 
    We design a rolling-horizon approach that optimizes the price menu whenever a new vehicle arrives, thereby tailoring incentives to real-time system conditions and individual user needs.  
    \item 
    Unlike prior menu-pricing studies that assume flat tariffs or ignore network constraints, our model simultaneously buys and sells energy at wholesale prices, supports bidirectional power flow, and enforces a feeder capacity constraint, capturing the full operational realities of a real parking lot.
    \item
    We embed the lower-level Karush–Kuhn–Tucker (KKT) conditions as indicator constraints to reformulate the bilevel model into a single mixed-integer linear program (MILP). In practice, this allows us to solve a full-day instance with 100 EVs in under 40 seconds and with 250 EVs in under 200 seconds, meaning each arriving EV can be priced in well under one second and optimized in real time.
    \item 
    Using 12 months of Australian Energy Market Operator (AEMO) price data, three tuned baseline schemes, and extensive sensitivity analysis (grid capacity, number of EVs, menu granularity), on average, the proposed mechanism increases parking lot operator profit by 27.01\%, reduces EV payments by 14.69\% and raises V2G contribution by 161.5\%.
\end{enumerate}

These results demonstrate that dynamic menu-based pricing offers a practical and economically efficient solution for parking lot operation in real-time electricity markets.

The remainder of the paper is structured as follows: Section~\ref{sec:related} reviews related work. Section~\ref{sec:model} describes the system model and Section~\ref{sec:Opt-problem} presents the problem formulation details and methodology.
Section~\ref{sec:case_study} discusses the case study results, followed by discussion and conclusion in Sections~\ref{sec:discussion} and~\ref{sec:conclusion}.

\section{Related Work}\label{sec:related}

Menu-based pricing has been widely studied as a mechanism to coordinate electric vehicle (EV) charging. The central idea is to offer drivers a set of price–service pairs (e.g., charging speed or deadline), thereby aligning user preferences with system-level objectives such as reducing costs, integrating renewables, or alleviating grid stress. Existing work varies in its assumptions about pricing dynamics, grid limits, and the availability of V2G discharging. We group related studies into four strands: (i) fixed menu designs under flat tariffs, (ii) online and adaptive pricing schemes, (iii) game-theoretic and bilevel formulations, and (iv) V2G integration and grid-aware models.

\emph{Fixed Menu Designs under Flat Tariffs:}
Bitar and Xu\cite{BitarXu2017Deadline} propose a deadline-based, charge-only menu that maximizes social welfare under a flat electricity tariff. Their scheme assumes one-way charging and excludes both V2G and price fluctuations.  
Wu et al.\cite{WuYucelZhou2021SmartCharging} also design a mechanism where each EV chooses a price–deadline pair, after which the operator schedules charging to reduce cost and emissions. While their results show significant efficiency gains, the model relies on day-ahead electricity prices, covers charging only, and does not consider feeder-capacity constraints.  
 Latinopoulos et al. \cite{Latinopoulos2021Optimal} take an empirical step by implementing menu pricing in a real parking lot. Every four hours, they publish a menu of 46 charging options or 60 V2G options, priced using genetic algorithms and particle swarm optimization. However, the menu remains fixed throughout the period, preventing adaptation to new arrivals or sudden wholesale price shifts.

\emph{Online and Adaptive Pricing Schemes}
Ghosh and Aggarwal \cite{GhoshAggarwal2018Menu} develop an online algorithm that recomputes a price list each time a new EV arrives. Drivers select energy–deadline bundles from this menu, which accounts for uncertain renewable output but assumes a flat electricity rate and excludes V2G.  
 Lu et al.\cite{Lu2023DeadlineMenu} extend deadline-sensitive menus by formulating the problem as a mixed-integer quadratic program calibrated with California EV data. While their model improves flexibility, it remains limited to charging, omits feeder-capacity constraints, and assumes exogenous tariffs. Mathioudaki et al.\cite{Mathioudaki2025OnlineMenu} propose an online tool that provides each arriving EV with a price–deadline offer. Their method guarantees incentive-compatible reporting but is restricted to charge-only operation, uses a fixed tariff, and excludes V2G and real-time price input.

\emph{Game-Theoretic and Bilevel Formulations}
Rasheed et al.\cite{Rasheed2024GameEVPricing} embed dynamic prices into a multi-leader, multi-follower Stackelberg game that coordinates EV routing and station selection. Although innovative, the model considers only one-way charging and ignores grid-capacity limits.  
Dupont et al.\cite{Dupont2024DRMenu} reformulate menu pricing in two mixed-integer linear models. One updates menus hourly to enable demand-response participation, while the other frames pricing as a bilevel game with incomplete information about drivers’ willingness to pay. Both assume unlimited chargers, flat energy tariffs, and charge-only operation, leaving out feeder constraints and V2G.
\emph{V2G Integration and Grid-Aware Models}
Ju and Moura\cite{JuMoura2023PricingScheme} introduce a three-option menu: fast charging, slower cheaper charging, and V2G with compensation for discharging. This simple design increases V2G participation under a time-of-use tariff, but the prices are fixed in advance, identical for all users, and do not account for feeder limits or real-time wholesale price changes.  
Latinopoulos et al.\cite{Latinopoulos2021Optimal}, as noted earlier, included V2G in their field experiment but relied on static four-hour menus.  
Overall, V2G integration has typically been explored only under static tariffs and without explicit consideration of grid constraints.
\paragraph*{Summary} In summary, most prior work relies on fixed menus or limited periodic updates, often under flat tariffs and with charging-only operation. Price updates are rarely performed in a rolling horizon, and almost no study combines real-time wholesale pricing, feeder-capacity limits, and full V2G participation in one model. Table~\ref{tab:comparison} highlights these gaps.  
Our work fills these gaps by introducing the dynamic menu pricing mechanism that jointly incorporates real-time wholesale prices, V2G, feeder constraints, and per-EV menus, and by providing a tractable MILP reformulation of the underlying bilevel problem.

\renewcommand{\arraystretch}{1.2}
\begin{table}[htbp]
\centering
\scriptsize
\caption{Feature support in prior studies and this work}\vspace{-2mm}
\label{tab:comparison}
\begin{tabularx}{\linewidth}{
  l
  *{6}{>{\centering\arraybackslash}X}
}
\toprule
\textbf{Study}
  & \makecell{\hspace*{-2mm}\textbf{On arrival}\\\hspace*{-2mm}\textbf{update}\strut}
  & \makecell{\textbf{V2G}\\\textbf{support}\strut}
  & \makecell{\hspace*{-2mm}\textbf{Real-time}\\\hspace*{-2mm}\textbf{prices}\strut}
  & \makecell{\textbf{Rolling}\\\textbf{horizon}\strut}
  & \makecell{\textbf{Per-EV}\\\textbf{menu}\strut}
  & \makecell{\hspace*{-2mm}\textbf{Network}\\\hspace*{-2mm}\textbf{Constraint}\phantom{}}\\
\midrule
Bitar and Xu\cite{BitarXu2017Deadline}
  & \xmark & \xmark & \xmark & \xmark & \xmark & \xmark \\
Ghosh and Aggarwal\cite{GhoshAggarwal2018Menu}
  & \xmark & \xmark & \xmark & \cmark & \cmark & \xmark \\
Latinopoulos et al.\cite{Latinopoulos2021Optimal}
  & \xmark & \cmark & \xmark & \xmark & \xmark & \cmark \\
Lu et al.\cite{Lu2023DeadlineMenu}
  & \xmark & \xmark & \xmark & \xmark & \cmark & \xmark \\
Rasheed et al.\cite{Rasheed2024GameEVPricing}
  & \xmark & \xmark & \xmark & \cmark & \cmark & \xmark \\
Dupont et al. \cite{Dupont2024DRMenu}
  & \xmark & \xmark & \xmark & \cmark & \xmark & \xmark \\
Ju and Moura\cite{JuMoura2023PricingScheme}
  & \xmark & \cmark & \xmark & \xmark & \xmark & \xmark \\
Wu et al.\cite{WuYucelZhou2021SmartCharging}
  & \xmark & \xmark & \xmark & \xmark & \cmark & \xmark \\
Mathioudaki et al.\cite{Mathioudaki2025OnlineMenu}
  & \xmark & \xmark & \xmark & \cmark & \cmark & \cmark \\
\midrule
\textbf{This work}
  & \cmark & \cmark & \cmark & \cmark & \cmark & \cmark \\
\bottomrule
\end{tabularx}
\end{table}

\section{System Model}
\label{sec:model}
We consider a parking lot operator that manages a set of bidirectional chargers connected to the distribution grid. 
The parking lot has unlimited parking/charging spaces, but is constrained by a feeder capacity $Y^{max}$ (kW). Each charger supports charging and discharging up to $P^{max}$ (kW).  
Time is discretized into slots $\mathcal{T} = \{1,\dots,T\}$ of equal duration $\tau \in \mathbb{R}^+$ hours. At each time slot $t \in \mathcal{T}$, the operator may trade electricity in the real-time market at buy and sell prices $\pi^{buy}_t$ and $\pi^{sell}_t$ (\$/kWh), respectively. We consider a set of EVs $\mathcal{N} = \{1,\dots,N\}$ arriving sequentially. 

\begin{figure*}
    \centering
    \includegraphics[width=0.9\textwidth]{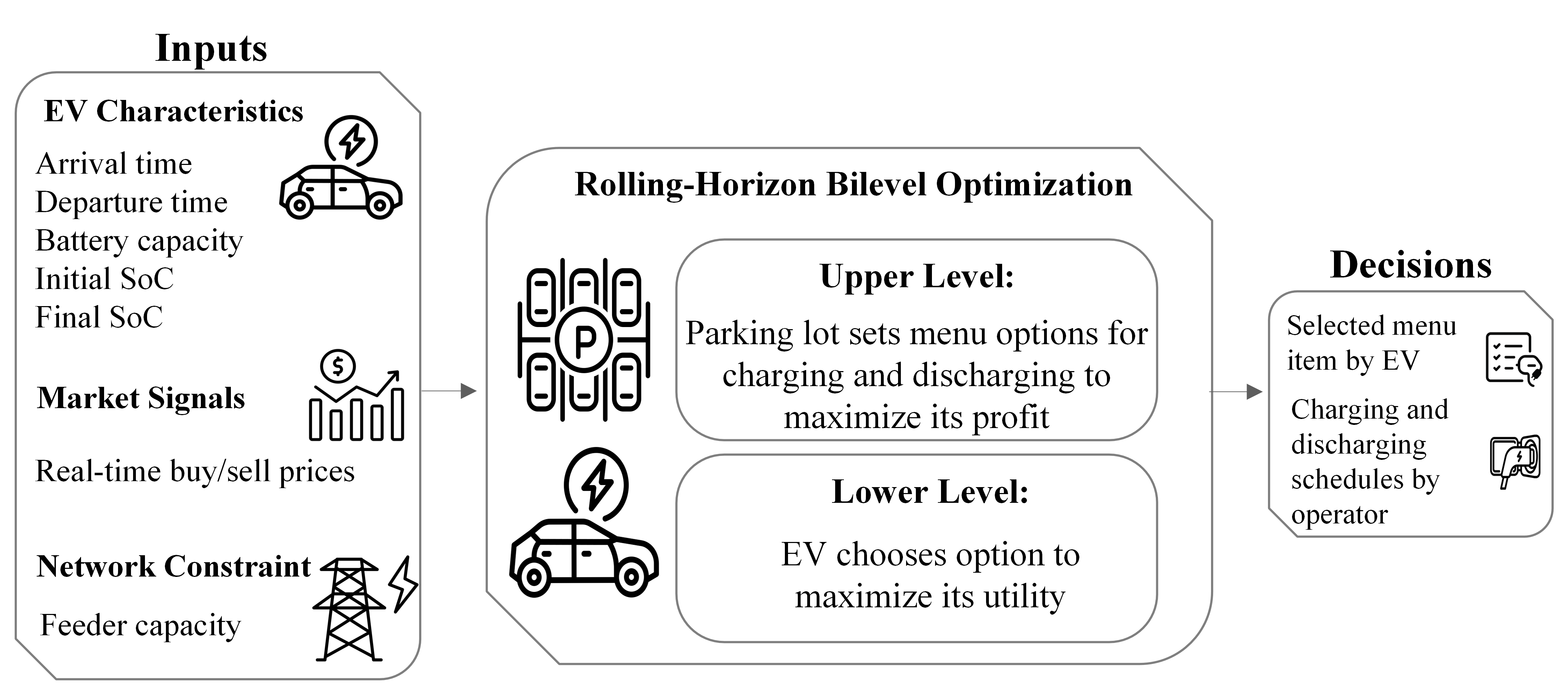}\vspace{-2mm}
    \caption{Overview of the proposed menu-based pricing framework}
    \label{fig:Overview}
\end{figure*}

When EV \(n\) arrives at time \(t_n^{\mathrm{arr}}\) with stated departure \(t_n^{\mathrm{dep}}\), it reports its battery capacity \(C_n^{\max}\) (kWh), its initial SoC, and final target SoC, respectively \(s^{{ini}}_n\) and \(s^{{fin}}_n\), scaled to \([0, 1]\). Upon this arrival, the operator prepares a  menu \(\mathcal{D}=\{d_1,\ldots,d_K\}\) of  discharge amounts options over \([t_n^{\mathrm{arr}},\,t_n^{\mathrm{dep}}]\). Here, each \(d\) (kWh) represents the total energy the EV would allow to be exported to the grid during its stay, and \(d=0\) corresponds to charge-only service. 

For each option \(d\in\mathcal{D}\), the operator computes a price \(\pi_{n,d}\) by solving a system optimization over the remaining horizon \([t_n^{\mathrm{arr}},\,T]\) that includes all active EVs \(\mathcal{E}_n\cup\{n\}\). Here, \(\mathcal{E}_n\) denotes the set of EVs admitted before \(n\) arrived; for each \(m\in\mathcal{E}_n\) the operator has a committed option \(d_m^{*}\) and a contracted price \(\pi_{m,d_m^{*}}\), which are not altered by later arrivals. In the optimization, the charging/discharging schedules of previously admitted EVs \(m\in\mathcal{E}_n\) may be adjusted within their technical bounds to restore feasibility and profitability under the option \(d\), while preserving those contracted prices. The optimization simultaneously enforces the feeder capacity \(Y^{\max}\), the per-charger limit \(P^{\max}\), and the SoC dynamics for every active EV.
\noindent In each optimization, the operator decides, for every time slot, how much power to import $P^{\mathrm{imp}}_{t,d}$ and how much power to export $P^{\mathrm{exp}}_{t,d}$, alongside EV charging/discharging. These quantities are chosen based on the buy/sell prices $\{\pi^{\mathrm{buy}}_t,\pi^{\mathrm{sell}}_t\}$ and must satisfy power balance and the feeder limit.

Once the menu \(\{\pi_{n,d}: d\in\mathcal{D}\}\) is computed, EV \(n\) responds by selecting the option that maximizes its own utility (or rejecting if all utilities are negative). If \(n\) accepts some \(d^{*}\), the corresponding price \(\pi_{n,d^{*}}\) and all schedules produced by that optimization are committed; the active set updates to \(\mathcal{E}_{n+1}=\mathcal{E}_n\cup\{n\}\). If \(n\) rejects, the schedules and prices for \(\mathcal{E}_n\) remain as previously committed, and the process advances to the next arrival. The mechanism is inherently online and rolling-horizon: it never assumes knowledge of future arrivals, it recomputes schedules only for the current active set at each arrival, and it preserves all price commitments already made to earlier EVs. At the arrival time $t_n^{\text{arr}}$ of each EV $n$, the optimization that generates the menu is built using the realized wholesale prices for the current slot $t_n^{\text{arr}}$ and historical or forecasted prices for all subsequent slots in the planning horizon. See Figure~\ref{fig:Overview} for an overview. This system model is the foundation for the bilevel formulation in Section~\ref{sec:Opt-problem}, where the operator’s profit-maximizing menu design (upper level) is coupled with the EV’s utility-maximizing option choice (lower level), and where the price for each option \(d\) is computed via optimization that incorporates energy prices, feeder limits, charger limits, and SoC feasibility. The resulting rolling-horizon implementation is summarized in Algorithm~\ref{alg:rolling-horizon-menu-pricing}.

\begin{algorithm}[!t]
\setstretch{0.95}
\caption{Rolling-Horizon MILP-Based Menu Pricing}
\label{alg:rolling-horizon-menu-pricing}
\begin{algorithmic}[1]
\Require Wholesale prices $\pi_t^{\text{buy}}$, $\pi_t^{\text{sell}}$; feeder and charger limits $Y^{\max}$, $P^{\max}$; sequence of arriving EVs $\mathcal{N}$
\Ensure Schedules $\{P^{\text{ch}}_{d^*}, P^{\text{dis}}_{d^*}\}$, offered prices $\{\pi_{n,d^*}\}$, and accepted EV set $\mathcal{E}_{n}$
\State $\mathcal{E}_0 \gets \emptyset$
\For{$n = 1$ to $|\mathcal{N}|$}
    \State $\mathcal{E}_n \gets \mathcal{E}_{n-1}$
    \State Read EV $n$'s profile: $(t_n^{\text{arr}}, t_n^{\text{dep}}, s_n^{\text{ini}}, s_n^{\text{fin}}, C_n^{\max})$
    \State \textit{// Compute baseline cost without new EV}
    \State $C^*(\mathcal{E}_{n-1}) \gets$ \Call{ComputeCost}{$\mathcal{E}_{n-1}$}
    \State \textit{// Upper level: Operator builds menu}
    \For{each option $d \in \mathcal{D}$}
        \State \textit{// Compute cost with new EV under option d}
        \State $C^*(\mathcal{E}_n \cup \{n\}) \gets$ \Call{ComputeCost}{$\mathcal{E}_n \cup \{n\}, d$}
        \State \textit{// Calculate marginal (incremental) cost}
        \State $mc_{n,d} \gets$ \Call{MarginalCost}{$\mathcal{E}_n \cup \{n\}, d$}
        \State \textit{// Set price with markup ($\beta_{n,d}$ is decision variable)}
        \State $\pi_{n,d} \gets mc_{n,d} + \beta_{n,d}$
        \State \textit{// Compute operator profit for this option}
        \State $\varphi_{n,d} \gets$ \Call{ComputeProfit}{$\pi_{n,d}, \mathcal{E}_n \cup \{n\}, d$}
        \State \textit{// Get charging and discharging schedules}
        \State $\{P^{\text{ch}}_d, P^{\text{dis}}_d\} \gets$ \Call{GetSchedules}{$\mathcal{E}_n \cup \{n\}, d$}
    \EndFor
    \State \textit{// Lower level: EV selects best option}
    \For{each option $d \in \mathcal{D}$}
     \State \textit{// Compute EV utility for this option}
     \State $U_{n,d} \gets$ \Call{ComputeEVUtility}{$n, d, \pi_{n,d}$}   
    \EndFor
    \State $d^* \gets \arg\max_{d \in \mathcal{D}} \{U_{n,d} : U_{n,d} \ge 0\}$
    \If{$d^*$ exists}
        \State \textit{// EV n accept the contract}
        \State $\mathcal{E}_n \gets \mathcal{E}_{n-1} \cup \{n\}$
    \Else
        \State \textit{// EV n reject the contract, keep previous set}
        \State $\mathcal{E}_n \gets \mathcal{E}_{n-1}$
    \EndIf
\EndFor
\State \Return $\mathcal{E}_{n}$, $\{\pi_{n,d^*}\}$, $\{P^{\text{ch}}_{d^*}, P^{\text{dis}}_{d^*}\}$
\end{algorithmic}
\end{algorithm}

\section{Bilevel Optimization Problem}
\label{sec:Opt-problem}

The coordination between the parking lot operator (leader) and the arriving EV $n$ (follower) is modeled as a bilevel optimization problem. The operator seeks to maximize their profit by strategically setting the menu of discharge options $\mathcal{D}$. Simultaneously, the EV selects the option from the menu that maximizes its personal utility. The optimal solution must respect the constraints of the system and the rational, utility-maximizing choice of the EV.
The problem is solved dynamically upon the arrival of EV $n$. The operator, as the leader, first determines the optimal menu prices ($\pi_{n,d}$) and the necessary charging/discharging schedules for a set of predefined discharge commitments $d \in \mathcal{D}$ (upper-level problem). The EV, observes the menu and selects the option $d^* \in \mathcal{D}$ that yields the highest non-negative utility (lower-level problem).

\subsection{Upper-Level Problem}
At the upper level, the parking lot operator (leader) maximizes the profit realized from the option chosen by the arriving EV $n$. For any candidate discharge option $d \in \mathcal{D}$ in the menu, the operator optimizes the system over the remaining horizon $[t_n^{\mathrm{arr}}, T]$. This optimization determines the EV charging and discharging schedules for all active EVs ($\mathcal{E}_n \cup \{n\}$), the grid exchanges, and the price $\pi_{n,d}$ offered to EV $n$ for that option.

The objective of the operator's problem is defined by two economic terms: (i) the price $\pi_{n,d}$ charged to EV $n$ for option $d$, and (ii) the wholesale settlement arising from the optimized grid imports and exports. The operator's profit $\varphi_{n,d}$ for a specific option $d$ is therefore:
\begin{equation}
\varphi_{n,d} = \pi_{n,d} - \tau \sum_{t=t_n^{\mathrm{arr}}}^{T} \left( \pi_t^{\mathrm{buy}} P^{\mathrm{imp}}_{t,d} - \pi_t^{\mathrm{sell}} P^{\mathrm{exp}}_{t,d} \right)
\end{equation}
where $\tau$ is the slot duration and $(\pi_t^{\mathrm{buy}}, \pi_t^{\mathrm{sell}})$ are the wholesale buy and sell prices, respectively. The schedules are decision variables that ensure feasibility and minimize the settlement term under the network and battery constraints.
Price $\pi_{n,d}$ charged to the EV is defined as a marginal-cost-based price: for each option $d$, the operator computes the incremental optimal settlement cost of admitting EV $n$ under option $d$ relative to the already-committed fleet $\mathcal{E}_n$ and sets the price $\pi_{n,d}$ based on a non-negative markup $\beta_{n,d}$:
\begin{equation}
\pi_{n,d} = \mathrm{mc}_{n,d} + \beta_{n,d}, \quad \text{with} \quad \beta_{n,d} \geq 0
\end{equation}
Here, $\beta_{n,d}$ is a decision variable chosen by the operator.The marginal cost $\mathrm{mc}_{n,d}$ is defined as the difference in the optimal system grid settlement costs:
\begin{equation}
\mathrm{mc}_{n,d} = \mathcal{C}^*(\mathcal{E}_n \cup \{n\}, d) - \mathcal{C}^*(\mathcal{E}_n)
\end{equation}
where $\mathcal{C}^*(\cdot)$ is the minimum grid settlement cost over the remaining horizon, formally defined as:
\begin{equation}
\mathcal{C}^*(\cdot) = \min \quad \tau \sum_{t=t_n^{\mathrm{arr}}}^{T} \left( \pi_t^{\mathrm{buy}} P^{\mathrm{imp}}_{t,d} - \pi_t^{\mathrm{sell}} P^{\mathrm{exp}}_{t,d} \right)
\end{equation}
This ensures that the price for option $d$ at least recovers the grid cost induced by the new EV. \\
We introduce a binary variable $z_{n,d} \in \{0, 1\}$ to indicate the EV's selection of option $d$. 
The upper-level objective is to maximize the expected realized profit, which is the sum of all option profits weighted by the selection variables:
\begin{equation}
\max \quad \sum_{d \in \mathcal{D}} z_{n,d} \varphi_{n,d}.
\end{equation}

We then introduce a scalar variable $\varphi^{\mathrm{s}}_n$ to represent the realized profit associated with the option chosen by EV $n$ and use indicator constraints to link it to the selected option:
\begin{align}
\label{eq:indicator_profit}
z_{n,d} = 1 \;\;\implies\;\; \varphi^{\mathrm{s}}_n = \varphi_{n,d},
\qquad \forall d \in \mathcal{D}.
\end{align}
If all options yield negative utility, the EV rejects the menu; in this case, all $z_{n,d} = 0$, we set $\varphi^{\mathrm{s}}_n = 0$, and the previously committed schedules and active EV set $\mathcal{E}_n$ remain unchanged.

Using $\varphi^{\mathrm{s}}_n$, the upper-level objective can equivalently be written as
\begin{equation}
\label{eq:obj2}
\max \;\; \varphi^{\mathrm{s}}_n.
\end{equation}
\
\subsubsection*{Operator Constraints}

The upper-level problem is subject to the following constraints.

The SoC of each EV evolves according to the net charged and discharged energy in each time slot:
\begin{equation}
\begin{aligned}
    s_{i,t,d} &= s_{i,t-1,d} 
    + \frac{\tau}{C_{i}^{\max}} 
    \left( 
        \eta_{\text{ch}}\, p^{\text{ch}}_{i,t,d} 
        - \frac{1}{\eta_{\text{dis}}}\, p^{\text{dis}}_{i,t,d}
    \right), \\
    &\forall i \in \mathcal{E}_n \cup \{n\},\ \forall d \in \mathcal{D},\ 
      \forall t \in \mathcal{T}.
\end{aligned}
\end{equation}
Here, $\eta_{\text{ch}}, \eta_{\text{dis}} \in (0,1]$ denote the charging and discharging efficiencies.

The SoC must start from the observed initial value and reach at least the required final value:
\begin{align}
    s_{i,t_i^{\text{arr}},d} &= s_i^{\text{ini}}, 
    && \forall i \in \mathcal{E}_n \cup \{n\},\ \forall d \in \mathcal{D}, \\
    s_{i,t_i^{\text{dep}},d} &\ge s_i^{\text{fin}}, 
    && \forall i \in \mathcal{E}_n \cup \{n\},\ \forall d \in \mathcal{D}.
\end{align}

For the arriving EV $n$, the total discharged energy under option $d$ must not exceed the corresponding commitment:
\begin{equation}
    \sum_{t=t_n^{\mathrm{arr}}}^{\mathcal{T}} \tau \, p_{n,t,d}^{\mathrm{dis}} 
    \;\leq\; d,
    \qquad \forall d \in \mathcal{D}.
\end{equation}

For each previously accepted EV $m \in \mathcal{E}_n$, the additional discharged energy over the remaining horizon is limited by its residual discharge budget:
\begin{equation}
\begin{aligned}
    \sum_{t = t_n^{\mathrm{arr}}}^{\mathcal{T}} 
        \tau \, p_{m,t,d}^{\mathrm{dis}}
    &\;\leq\; \bar d^{*}_m, \\
    &\forall m \in \mathcal{E}_n,\ \forall d \in \mathcal{D},
\end{aligned}
\end{equation}
where $\bar d^{*}_m \ge 0$ is the residual discharge budget of EV $m$ at the current decision step.

Charging and discharging powers for each EV are limited by the charger capacity:
\begin{align}
    0 \;\leq\; p^{\text{ch}}_{i,t,d},\ p^{\text{dis}}_{i,t,d}    &\leq P^{\max}, 
    && \forall i \in \mathcal{E}_n \cup \{n\},\ \forall t \in \mathcal{T},\ \forall d \in \mathcal{D}.
\end{align}

Each EV can only exchange power with the charger while it is physically present at the parking lot:
\begin{align}
    p^{\text{ch}}_{i,t,d},\ p^{\text{dis}}_{i,t,d} &= 0,
    && \forall i \in \mathcal{E}_n \cup \{n\},\ \forall d \in \mathcal{D},\
       \forall t \in \mathcal{T} \setminus [t_i^{\text{arr}}, t_i^{\text{dep}}].
\end{align}

At each time slot and for each option, the total power entering the parking lot must equal the total power leaving it:
\begin{equation}
\begin{aligned}
    \tau P^{\mathrm{imp}}_{t,d} 
    &+ \sum_{i \in \mathcal{E}_n \cup \{n\}} \tau p^{\text{dis}}_{i,t,d} 
    = \tau P^{\mathrm{exp}}_{t,d} 
    + \sum_{i \in \mathcal{E}_n \cup \{n\}} \tau p^{\text{ch}}_{i,t,d}, \\
    &\forall t \in \mathcal{T},\ \forall d \in \mathcal{D}.
\end{aligned}
\end{equation}
Moreover, at each time slot, grid import and export are mutually exclusive, and each EV can either charge or discharge, but not both simultaneously.

The power exchanged with the grid is constrained by the feeder capacity:
\begin{align}
    0 \;\leq\; P^{\mathrm{imp}}_{t,d},\ P^{\mathrm{exp}}_{t,d} &\leq Y^{\max}, 
    && \forall t \in \mathcal{T},\ \forall d \in \mathcal{D}.
\end{align}

Finally, the arriving EV can select at most one option from the menu:
\begin{equation}
    \sum_{d \in \mathcal{D}} z_{n,d} \;\leq\; 1.
\end{equation}

\subsection{Lower-Level Problem}
The arriving EV $n$ acts as the Follower, making a rational choice by selecting the option $d^*$ from the menu $\mathcal{D}$ that maximizes its personal utility $U_{n,d}$. This selection represents the Lower-Level Problem, which must be satisfied as a constraint in the overall bilevel formulation. The EV only selects an option if the resultant utility is non-negative.

The EV's utility $U_{n,d}$ is defined as the total perceived value derived from the required charge minus the total cost of the service. The cost comprises the price $\pi_{n,d}$ and the perceived cost of battery degradation from discharging:
\begin{equation}
U_{n,d} = \alpha_{n} \cdot (s^{\mathrm{fin}}_{n}-s^{\mathrm{ini}}_{n}) \cdot C^{\max}_{n} - \pi_{n,d} - \gamma_{n} d
\end{equation}
Here, $\alpha_{n}$ ($\$/\mathrm{kWh}$) is the EV’s energy valuation, and $\pi_{n,d}$ is the price determined by the operator (Leader). The term $\gamma_{n} d$ is the perceived cost of discharging $d$ kWh, where $\gamma_{n}$ is the per-unit battery degradation cost ($\$/\mathrm{kWh}$), which is estimated using the model of \cite{Qi2025},  expressed as:
\vspace{-1mm}\begin{align}
\gamma_{n} = \frac{R}{L \cdot e}.
\label{eq:deg}
\end{align}
In this estimation, $R$ is the battery replacement cost, $L$ is the lifetime energy throughput, and $e$ is the square root of the roundtrip efficiency. Multiplying $\gamma_{n}$ by $d$ gives the total degradation cost when the EV discharges $d$~kWh. 
 
 The EV's decision rule, the Lower Level Problem, is to choose the option $d^*$ that yields the maximum non-negative utility:
\vspace{-2mm}\begin{equation}
z_{n,d} \in \arg\max_{d \in \mathcal{D}} \{ U_{n,d} \mid U_{n,d}\ge 0 \}
\end{equation}
The binary variable $z_{n,d^*}=1$ if option $d^*$ is chosen, and $z_{n,d}=0$ otherwise.
\subsection{Lower Problem Reformulation}
To efficiently solve this bilevel optimization problem, we reformulate it into a single-level Mixed-Integer Linear Programming (MILP) problem~\cite{Audet1997Links}. This is achieved by incorporating the Karush-Kuhn-Tucker (KKT) conditions of the lower-level problem and translating them into linear constraints using indicator variables.

The lower-level problem involves the EV selecting a discharge option \( d \in \mathcal{D} \) to maximize its utility. Although the variables $z_{n,d}$ represent binary choices, we relax them to $z_{n,d}\geq 0$ with $\sum_{d} z_{n,d}\leq 1$ so that the lower-level problem is a linear program and strong duality applies. This relaxation is exact: in 
generic cases the solution assigns one $z_{n,d}=1$ and all others zero; in degenerate cases multiple options may tie, but at least one binary optimal solution always exists. The resulting primal formulation is:
\begin{align}
\max \quad & \sum_{d \in \mathcal{D}} z_{n, d} U_{n,d}, \\
\text{s.t.} \quad 
&  \sum_{d \in \mathcal{D}} z_{n,d} \leq 1,\\
&\quad z_{n,d} \geq 0 \quad \forall d \in \mathcal{D}.
\label{eq:none_or_one}
\end{align}
From this primal problem, we derive the following dual:
\begin{align}
\min \quad & \lambda \\
\text{s.t.} \quad 
& \lambda \geq U_{n,d}, \quad \forall d \in \mathcal{D}, \\
& \lambda \geq 0. \label{eq:dual}
\end{align}
The KKT optimality conditions include complementary slackness:
\begin{align}
z_{n,d} (\lambda - U_{n,d}) = 0, \quad \lambda \left( \sum_{d \in \mathcal{D}} z_{n,d} - 1 \right) = 0, \quad \forall d \in \mathcal{D}. \label{eq:complementarity}
\end{align}
To reformulate the nonlinear product terms in the KKT conditions, we use indicator constraints. For each discharge option \( d \in \mathcal{D} \), we enforce:
\begin{align}
z_{n, d} = 1 \implies \lambda = U_{n,d}, \sum_{d \in \mathcal{D}} z_{n,d} = 0 \implies \lambda = 0. \label{eq:indicator}
\end{align}
These constraints ensure equivalence to the KKT conditions and enable the bilevel problem to be reformulated as a single-level MILP.

\section{Experimental Design}
\label{sec:case_study}
\subsection{Simulation Setup}
\label{sec:setup}
This experiment is designed to evaluate the proposed menu-based pricing algorithm for a simulated parking lot.

\begin{figure}
\centering
\includegraphics[width=0.5\columnwidth]{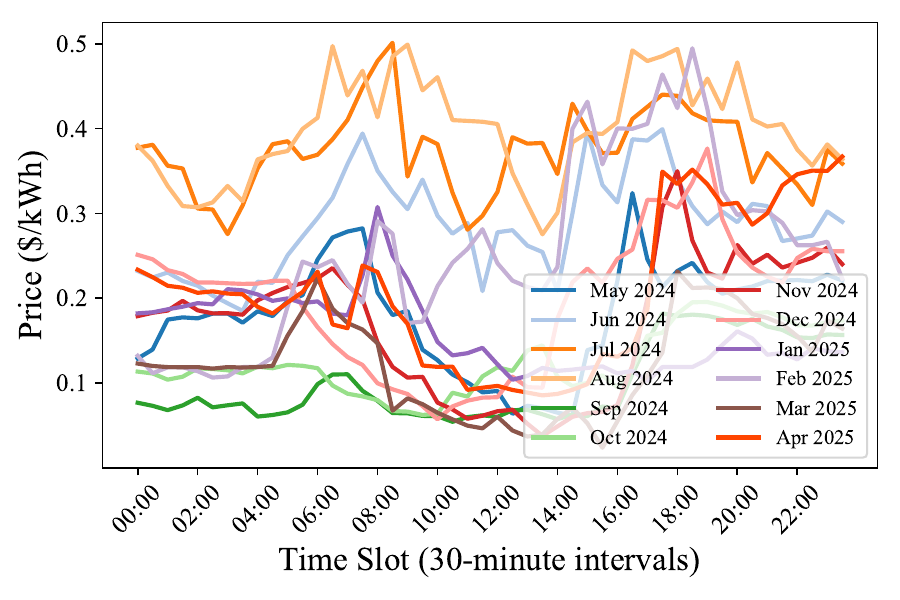} \vspace{-2mm}
\caption{Price profile of AEMO for Monday of each month.}
\label{fig:12monthprice}
\end{figure}

To simulate the parking lot, the day is divided into \(T = 48\) half-hour slots, providing a 24-hour planning horizon. Each simulated day contains \(N = 100\) EVs, each with a \(60~\mathrm{kWh}\) battery capacity of mid-range EVs such as the Nissan Leaf.  Initial and final SoC are sampled from truncated normal distributions~\cite{Sagaria2025V2GCompensation} with means 0.30 and 0.80 (standard deviation 0.10), reflecting the mix of commuters and shoppers in the urban area.  EV arrivals are modeled using a time-varying Poisson process to reflect realistic patterns of urban parking demand and parking durations for each EV are sampled from a uniform distribution between 2 and 6 hours.
Bidirectional chargers are rated at \(P^{\max}=60~\mathrm{kW}\), while the
feeder capacity is limited to \(Y^{\max}=600~\mathrm{kW}\).
In our simulations, grid settlements are assumed to follow wholesale prices published by the Australian Energy Market Operator (AEMO)\footnote{\url{https://aemo.com.au}}, with these prices accessed through Amber Electric\footnote{\url{https://www.amber.com.au}}, a retailer that passes real-time wholesale signals to end users.
The parking lot operator pays $\pi_t^{\text{buy}} = \lambda_t + \delta$ when importing and receives $\pi_t^{\text{sell}} = \lambda_t$ when exporting, 
where $\lambda_t$ is the real-time wholesale price. The extra term $\delta$ represents retailer and network costs, 
applied only to imports, while exports are credited at the wholesale price. To reduce computational complexity while still capturing seasonal variation, we select the first Monday of each month from our 12-month dataset (May 2024–April 2025) as representative daily price profiles for simulation (Figure~\ref{fig:12monthprice}).
In all simulations, we assume that the operator has access to a forecast of wholesale prices for the remainder of the day. This is consistent with common practice in electricity markets, where operators rely on publicly available price forecasts or retailer-provided predictions. These forecasted prices are treated as the operator’s best estimate of future wholesale conditions and are used when computing menu prices and evaluating each EV’s charging/discharging schedule over the remaining horizon. This perfect-forecast assumption allows us to isolate the performance of the pricing mechanism itself, without confounding effects from price-prediction uncertainty; the impact of forecast errors is examined separately in the sensitivity analysis of Section~\ref{sec:result}.

In practice, an EV user’s valuation of stored energy can be inferred from historical charging behavior (e.g., home charging costs, typical willingness to pay, and preferences observed in smart-meter data). Once deployed, such valuations can be personalized using real user profiles. However, because our study focuses on an Australian urban setting and no individual behavioral data are available, we adopt a representative value based on typical residential electricity tariffs. Most Australian households face retail prices around $0.25–0.35/kWh$, so we set the energy valuation parameter to a conservative and widely used benchmark of $\alpha_n = 0.30\,$\$/kWh. Under this assumption, any menu offering a charging price above this threshold yields negative utility and is therefore rejected by the user.
For battery degradation, we rely on an engineering-based calculation rather than an assumed behavioral parameter. Using Equation (19), we take a 60\,kWh Nissan Leaf battery with approximately 1,500 full cycles to an 80\% end-of-life threshold and a round-trip efficiency of 90\%. This gives a lifetime throughput of $L = 60 \times 1500 = 90{,}000$\,kWh and an efficiency factor of $e = \sqrt{0.90} \approx 0.95$. With a replacement cost of $R = 12{,}000$\,AUD, the implied degradation cost is $\gamma_n \approx 0.14\,$\$/kWh. 
Although our simulations use common values for all EVs, the framework can seamlessly incorporate heterogeneous $\alpha_n$ and $\gamma_n$ using either real historical charging patterns from participating EV owners, or manufacturer-level degradation data. These inputs can be updated dynamically once real-world pilot data become available.

\subsection{Baseline Methods for Comparison}
\label{sec:baselines}
\begin{figure}
\centering
\includegraphics[width=0.5\columnwidth]{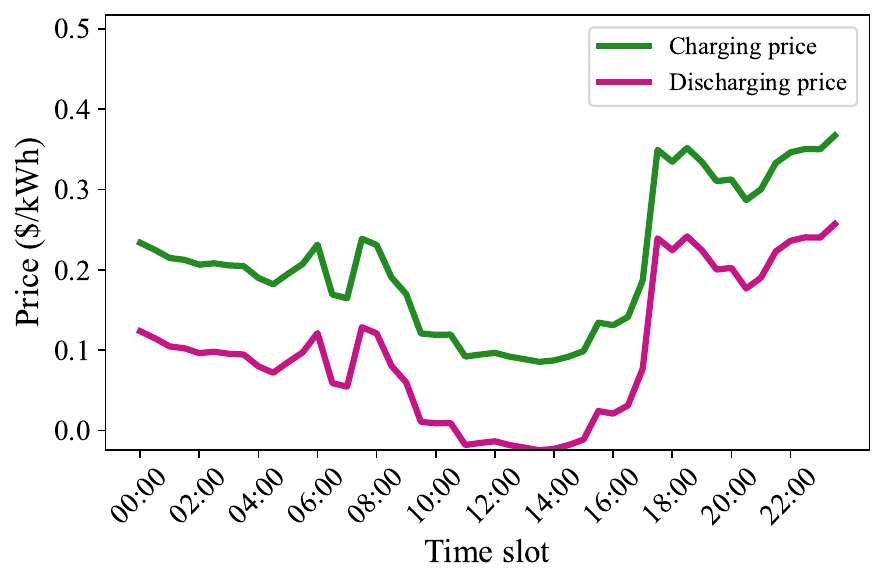} \vspace{-2mm}
\caption{Price profile for trading energy with grid; first Monday of April 2025.}
\label{fig:amberPricing}
\end{figure}

To evaluate the performance of our proposed menu-based pricing mechanism, we compare it against three simplified baseline strategies. In all cases, the parking lot trades energy with the grid through a retailer (e.g., Amber Electric), using the same real-time prices for both buying and selling energy, as shown in Figure~\ref{fig:amberPricing}.

\begin{figure}
\centering
\includegraphics[width=0.8\textwidth]{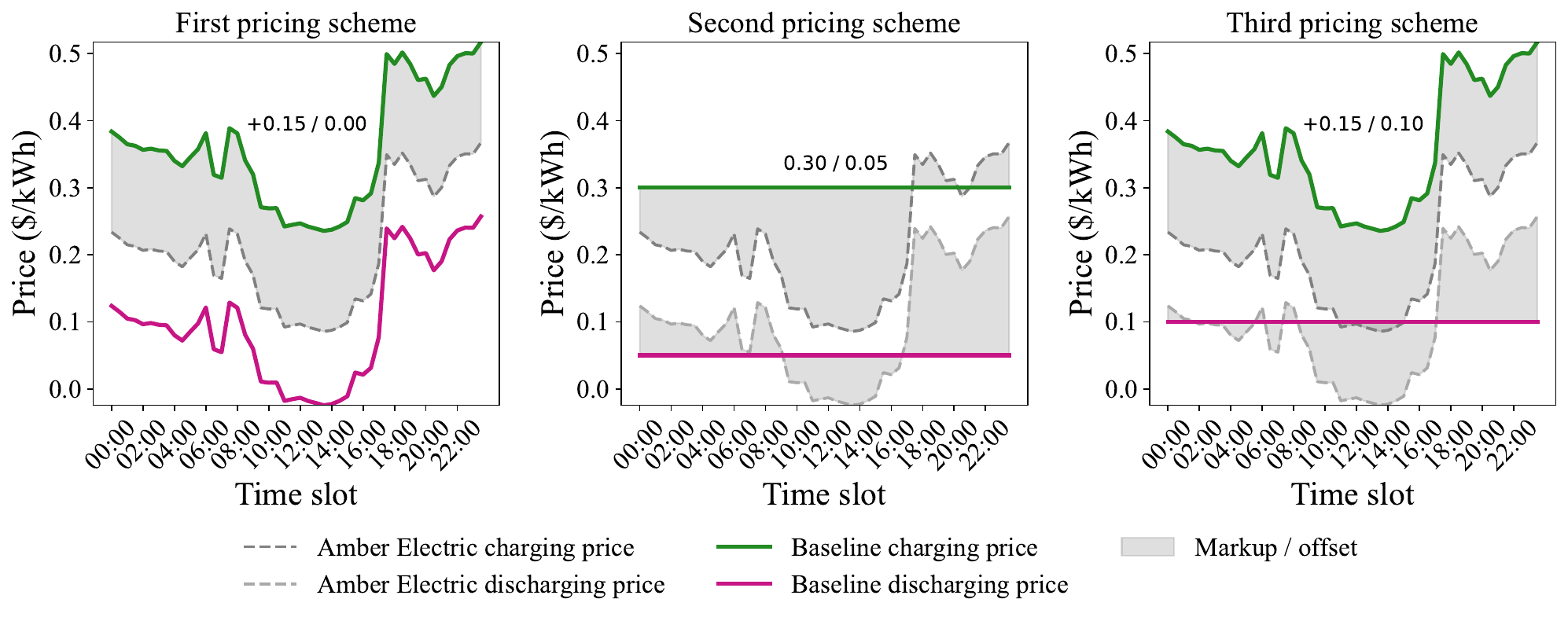}\vspace{-2mm}
\caption{Charging and discharging price profiles offered to EV users under the three baseline schemes.}
\label{fig:diffPricingProfile}
\end{figure}

Therefore, the grid-side prices are consistent across all methods. The key differences lie in how charging and discharging prices are offered to the EV users, as shown in Figure~\ref{fig:diffPricingProfile}. This ensures that any performance differences arise only from the pricing strategy offered to EV users, rather than from differences in grid-side transactions.

\paragraph{\textbf{Baseline 1: Adjusted Real-time Pricing}}
Charging prices offered to EVs are based on the real-time price plus a fixed markup to account for network and retail costs. Discharging prices offered to EVs are set to the price minus an offset.

\paragraph{\textbf{Baseline 2: Fixed Flat Tariff}}\
Both charging and discharging prices offered to EVs are set to fixed, time-invariant values throughout the day. This reflects a basic flat-rate scheme, independent of wholesale market conditions.

\paragraph{\textbf{Baseline 3: Hybrid Pricing}}\
Charging prices follow the same adjusted formula as in Baseline 1, while discharging prices are set to the fixed value. This creates a hybrid approach where charging reflects market variation but discharging does not.

For each baseline method, we optimize the price markup parameters to ensure a fair and meaningful comparison with our proposed menu-based mechanism. Specifically, for each baseline, we conduct a grid search by varying both the charging and discharging markups from 0 to 0.30 in steps of 0.05. For each combination of markup values, we simulate the parking lot operation and record the resulting profit and user payment. The markup upper bound of 0.30 \$/kWh is chosen to ensure that the final price offered to each EV remains within the user's maximum acceptable valuation defined in Section~\ref{sec:setup}. Prices above this level would provide no incentive for EVs to accept charging or discharging options.

Because of this upper bound, we observed that all markup combinations resulting in final prices above 0.30 \$/kWh  were consistently rejected by users, yielding zero accepted offers and thus zero parking lot profit. Therefore, these cases are excluded from the experiments, as they provide no meaningful contribution to the analysis. For all comparisons, we report the baseline results using the markup combination that yields the highest parking lot profit.

\section{Results}
\label{sec:result}
\subsection{Overall Performance}
\begin{figure}[!ht]
\centering
\includegraphics[width=0.9\textwidth]{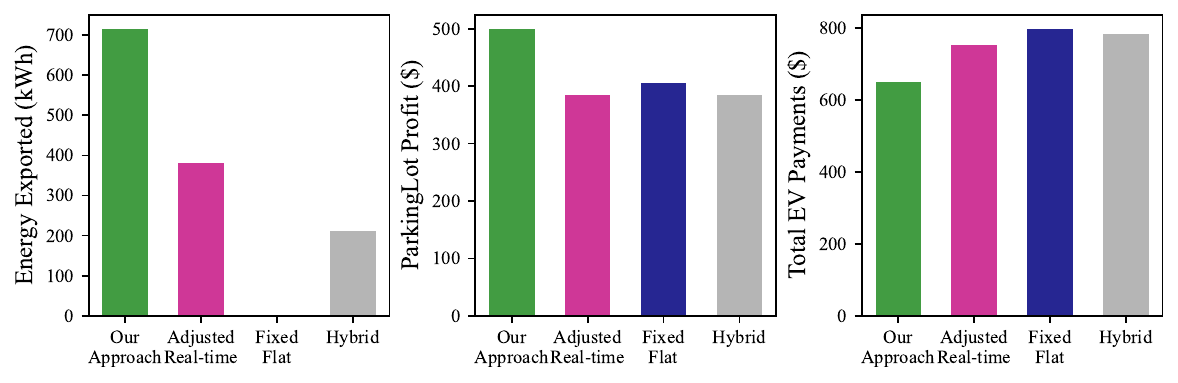} \vspace{-3mm}
\caption{Overall performance of our menu-based mechanism versus three baselines.}
\label{fig:totalcompare}
\end{figure}

Figure~\ref{fig:totalcompare} shows that our proposed menu-based pricing mechanism achieves better outcomes than all three baselines. The total energy exported to the grid is noticeably higher, indicating that EVs are more willing to participate in discharging under our dynamic pricing scheme. This increased V2G activity not only supports the grid but also leads to higher profit for the parking lot operator. At the same time, the total cost paid by EV users is lower, as they receive tailored incentives that reduce their payment. Overall, the figure highlights that our mechanism benefits both the operator and the users while improving grid support. 

\begin{figure}[!ht]
\centering
\includegraphics[width=0.4\columnwidth]{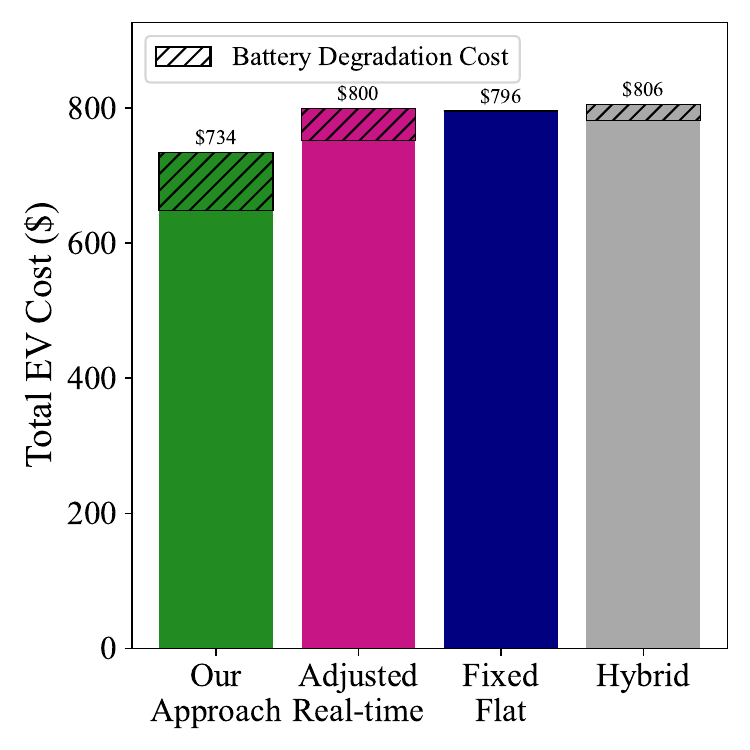}\vspace{-2mm}
\caption{EV Cost Breakdown: Energy payment and battery degradation cost for our mechanism and baselines.}
\label{fig:BatteryCost}
\end{figure}

Figure~\ref{fig:BatteryCost} breaks down the cost that each EV
pays into two parts: the energy cost and the extra amount we add for
battery degradation.  Even after this degradation cost is included, the menu-based scheme still gives the lowest overall cost, while the baselines stay higher.  This means our pricing keeps its advantage even when the real cost of battery degradation is counted.

\begin{figure*}[t]
\centering
\includegraphics[width=0.98\textwidth]{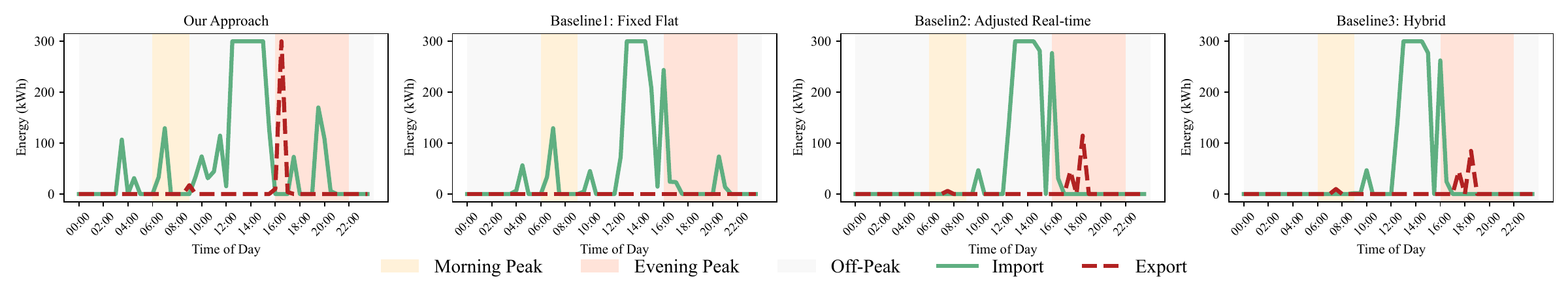}
\caption{Grid import and export profiles for our mechanism and baselines across a 24-hour day.}
\label{fig:ExIm}
\end{figure*}

Figure~\ref{fig:ExIm} illustrates the import and export energy profiles for our proposed mechanism compared to three baseline methods across a 24-hour period, segmented into morning peak (6-8am), evening peak (4-10pm), and off-peak times. During peak periods, our mechanism distinctly shows higher energy exports and lower energy imports compared to the baseline scenarios. Conversely, during off-peak periods, the pattern reverses, demonstrating lower exports and higher imports, capitalizing on lower electricity prices to economically charge the vehicles. This dynamic adjustment not only provides substantial economic benefits to both the parking lot operator and EV users but also enhances grid support through increased energy exports precisely when needed most.

\subsection{Sensitivity Analysis}
In order to assess the robustness of our proposed mechanism, we conduct a sensitivity analysis. The goal of this experiment is not only to test whether the mechanism performs well under a fixed set of assumptions, but also to understand how key parameters influence the outcomes of interest. By varying design choices, we can observe how these changes affect metrics like operator profit, user payments, and V2G participation.
To show what V2G adds, we also test a charge-only version of our method. It uses the same data, constraints, prices, and rolling-horizon setup, but discharging are turned off and the menu is $\mathcal{D}=\{0\}$. We label this variant ``Charge-only'' in all figures and tables.
\subsubsection{\textbf{Effect of Menu Options on our Method}}
\begin{figure}[!ht]
\centering
\includegraphics[width=0.6
\columnwidth,trim={0 40 0 0},clip]{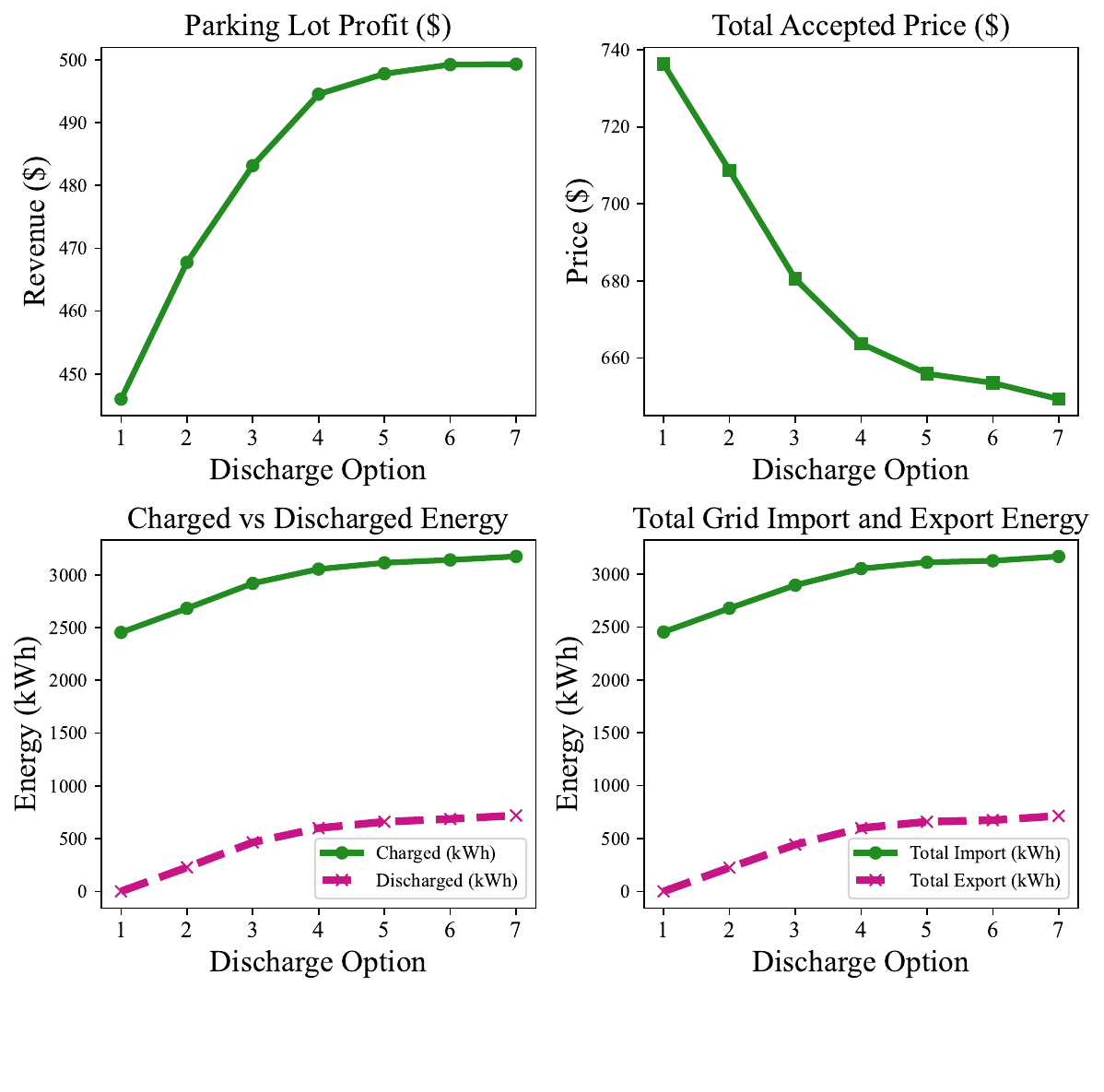} \vspace{-2mm}
\caption{Effect of menu size: profit, user payments, and discharged energy versus number of discharge options.}
\label{fig:optionImpact}
\end{figure}

\
Figure~\ref{fig:optionImpact} and Table~\ref{tab:optionImpactTable} show that when the menu expands from a minimal set to about 6–8 flexible discharge options, parking lot profit, user payments, and V2G participation all rise rapidly and then plateau. A richer menu lets the operator better match EV preferences and increase V2G participation, but adding more than roughly eight options brings little extra value. Thus, a small, well-chosen set of discharge levels captures most of the benefit while keeping the system simple.

\begin{table}[H]
    \centering
    \caption{Different menu
    option sets.}\vspace{-2mm}
    \scriptsize
    \begin{tabular}{cc}\toprule
         Menu & Discharge Options \\ \midrule
         1 & 0 \\
         2 & 0, 5 \\
         3 & 0, 5, 10 \\
         4 & 0, 5, 10, 15 \\
         5 & 0, 5, 10, 15, 20, 25, 30 \\
         6 & 0, 5, 10, 15, 20, 25, 30, 35, 40 \\
         7 & 0, 5, 10, 15, 20, 25, 30, 35, 40, 45, 50\\ \bottomrule
    \end{tabular}
    \label{tab:optionImpactTable}
\end{table}

\subsubsection{\textbf{Impact of Grid Capacity Limit}}
\begin{figure}[!ht]
\centering
\includegraphics[width=0.65\columnwidth]{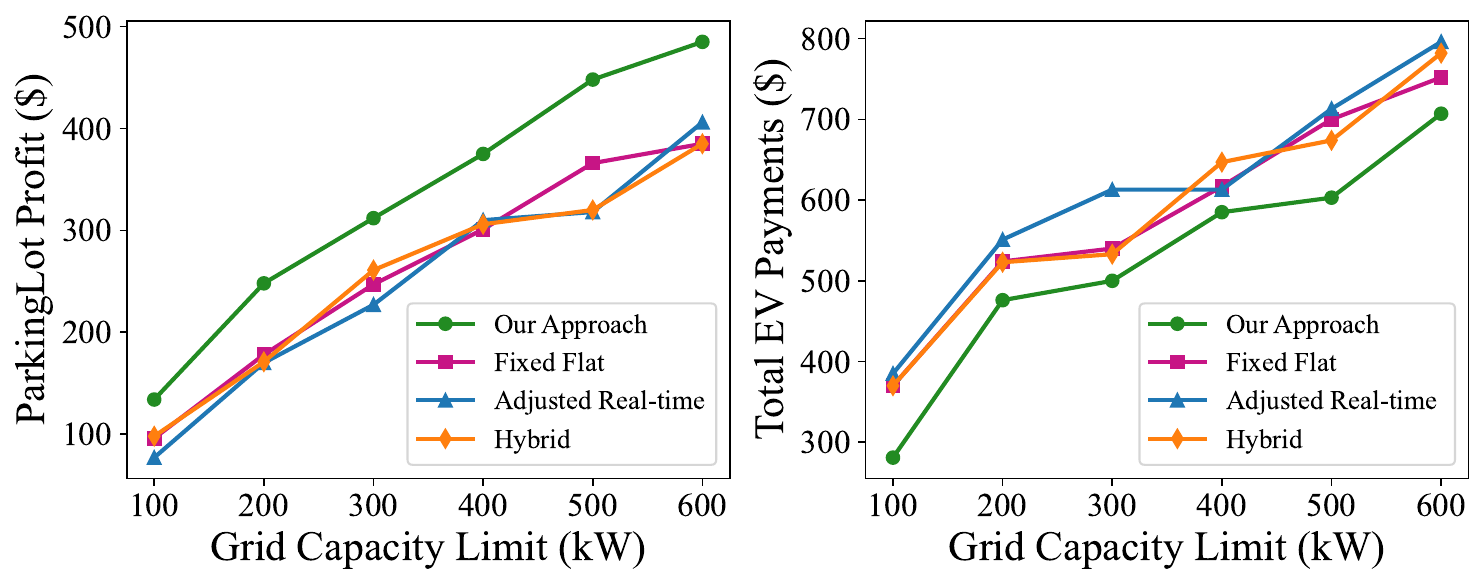} \vspace{-2mm}
\caption{Sensitivity to feeder-capacity limit, comparing our mechanism with all three baseline schemes.}
\label{fig:GridImpact}
\end{figure}

\
Figure~\ref{fig:GridImpact} explores how profit and user payments respond to varying feeder capacity limits. Our menu-based mechanism consistently yields higher parking lot profit and lower user costs, even as grid constraints become tighter. 

\subsubsection{\textbf{Effect of Number of EVs}}
\begin{figure}[!ht]
\centering
\includegraphics[width=0.65\columnwidth]{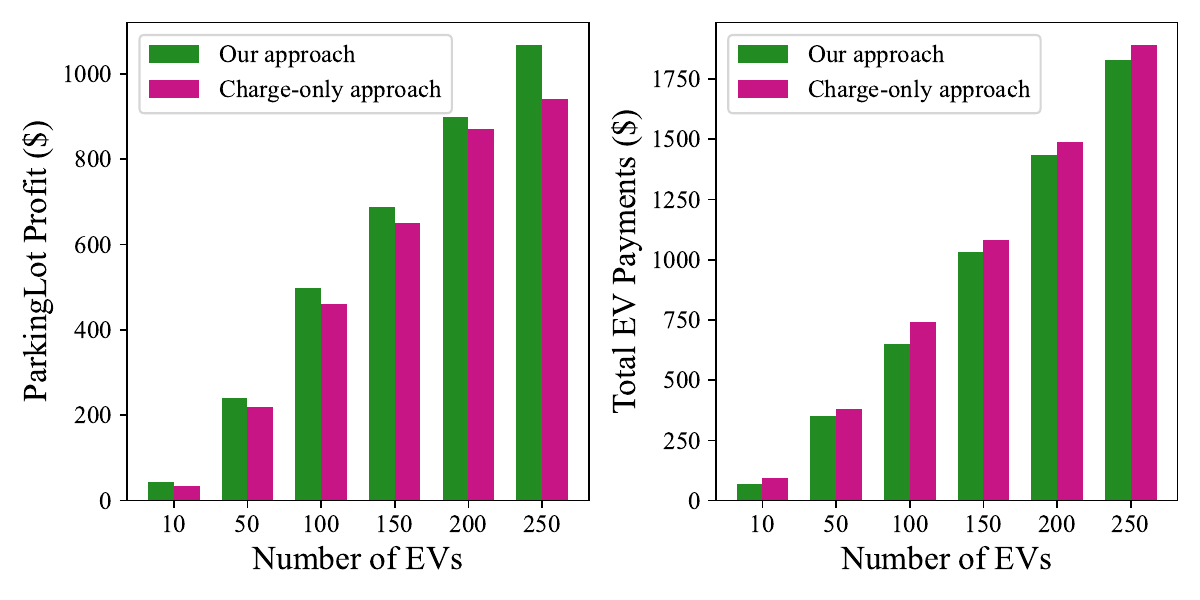} \vspace{-2mm}
\caption{Sensitivity to the number of participating EVs, comparing our mechanism with a charge-only approach.}
\label{fig:EVImpact}
\end{figure}

\
As the number of participating EVs increases (Figure~\ref{fig:EVImpact}), our mechanism remains effective, delivering superior outcomes in both profit and energy export compared to all baselines. This demonstrates that benefits persist in both small and large parking lots.

Moreover, the rolling optimization for each EV executes in under one second. In our experiments, all optimization problems were solved using Gurobi~12 on a MacBook Pro with an Apple M3 Pro chip and 18~GB memory. Under this setup, a full-day instance with 100 EVs was solved in under 40 seconds, while an instance with 250 EVs was solved in under 200 seconds. This sub-second processing when an EV arrives highlights the potential for real-world deployment by parking lot operators without the need for specialized infrastructure.

\subsubsection{\textbf{Monthly and Seasonal Variation}}
\begin{figure}[!ht]
\centering
\includegraphics[width=0.65\columnwidth]{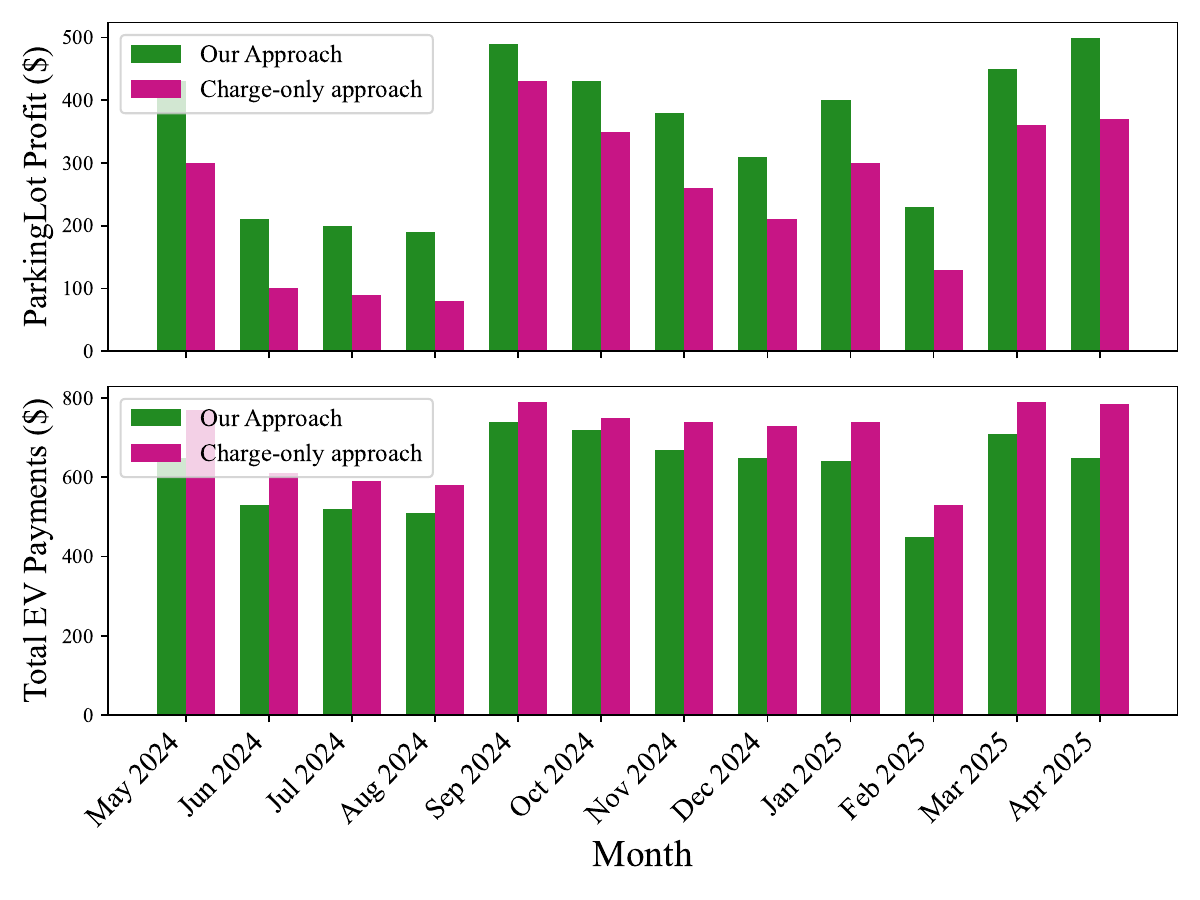} \vspace{-2mm}
\caption{Sensitivity to monthly wholesale-price variation (first Monday of each month), comparing our approach with a charge-only approach.}
\label{fig:MonthImpact}
\end{figure}

\
To assess resilience to seasonal effects, Figure~\ref{fig:MonthImpact} shows performance over the first Monday of each month. The menu-based strategy maintains its advantages year-round, consistently achieving higher parking lot operator profit and user cost savings, regardless of underlying fluctuations in market prices.
\subsubsection{\textbf{Robustness and Sensitivity to Wholesale Price Uncertainty}}
\begin{figure}[t]
\centering

\begin{subfigure}[b]{0.49\columnwidth}
    \centering
    \includegraphics[width=\linewidth]{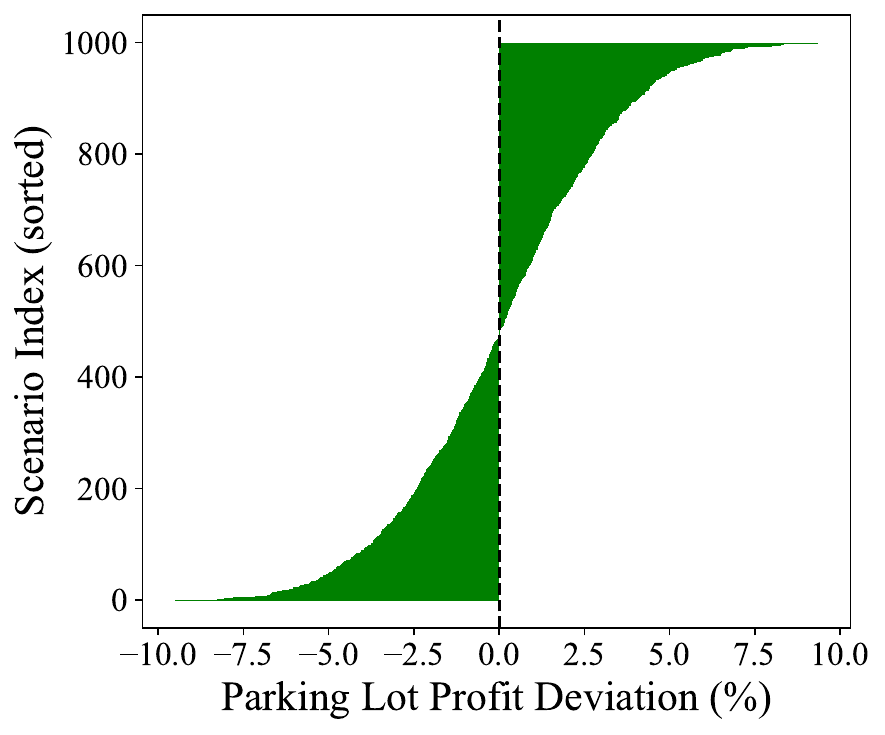}
    \caption{}  
    \label{fig:error}
\end{subfigure}
\hfill
\begin{subfigure}[b]{0.49\columnwidth}
    \centering
    \includegraphics[width=\linewidth]{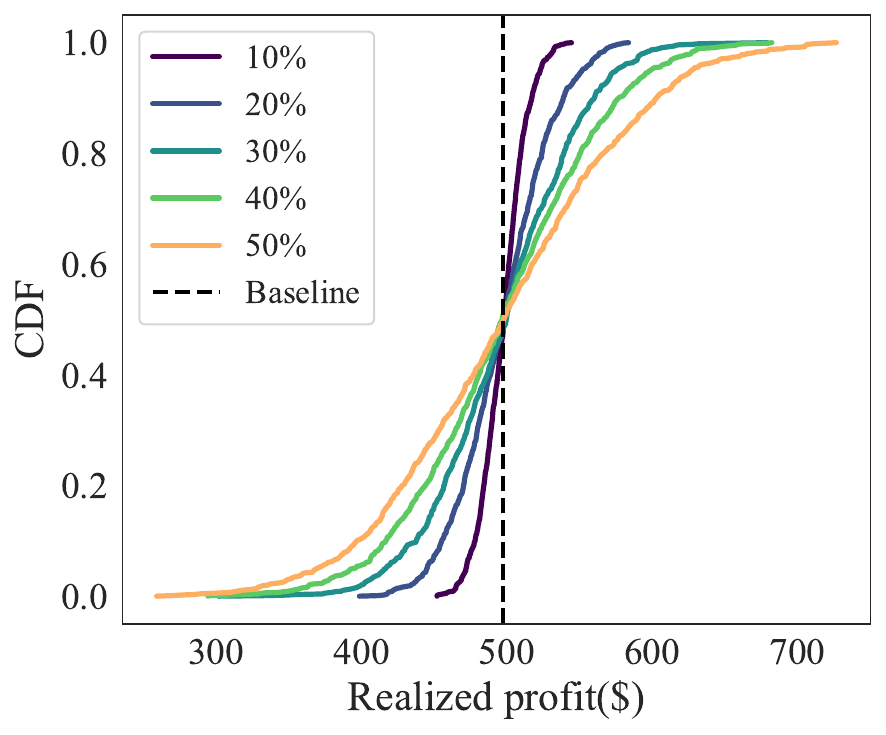}
    \caption{} 
    \label{fig:cdf}
\end{subfigure}
\vspace{-2mm}
\caption{Sensitivity to wholesale-price uncertainty: (a) deviation under 10\% noise; (b) CDF under 10--50\% noise.}
\label{fig:combined}
\end{figure}

In the baseline experiments, we assume that the operator has perfect knowledge (or highly accurate forecasting) of the wholesale price trajectory for the simulated day. This assumption is common in mechanism-design and EV-charging studies, as it isolates the behavior of the pricing mechanism without the influence of forecasting uncertainty. In practice, however, the operator schedules charging and discharging using a forecast, and the realized operating cost may differ if wholesale prices deviate unexpectedly. To quantify the effect of such deviations, we perform a Monte Carlo sensitivity analysis in which the operator’s charging and discharging schedules (optimized using the AEMO price profile) are held fixed, while the wholesale prices are perturbed ex post. This targets the question: \emph{if the operator schedules optimally based on a forecast, how much profit volatility does price uncertainty introduce at realization time?}\\
\vspace{-0.1 mm} 
We generate $1000$ perturbed price scenarios by adding Gaussian noise with standard deviation $10\%$ to each time-slot price. For each noisy scenario, we recompute only the operator profit. Across these $1000$ scenarios, the mean absolute percentage deviation is $2.42\%$, indicating modest variability under price uncertainty. The probability that the profit decreases by more than $5\%$ across all scenarios is just $4.9\%$. Figure~\ref{fig:error} shows the distribution of percentage deviations across all noisy-price scenarios. 
Overall, the proposed menu-based pricing mechanism exhibits low sensitivity to wholesale price uncertainty, demonstrating that it maintains reliable operational and economic performance under realistic forecasting errors.\\
\vspace{-0.1 mm} 
Figure ~\ref{fig:cdf} further shows the Cumulative Distribution Function (CDF) of realized operator profit under increasing wholesale-price forecast errors. As the noise level rises from 10\% to 50\%, the distribution gradually widens. However, the median profit remains remarkably stable across all uncertainty levels (496--499\,\$), staying almost identical to the baseline value of 497.29\,\$. These results indicate that although higher price uncertainty increases variability, the mechanism preserves its central profit performance and avoids severe downside outcomes, demonstrating robustness to substantial forecast errors.\\
\vspace{-0.1 mm} 
Table~\ref{tab:summary_improvement} sums up the core results.
Compared with the baselines, our menu with V2G exports much more energy, earns higher profit for the parking lot, and lowers the total amount paid by EV owners.
In this table, values indicate desirable improvements, higher exported energy and parking-lot profit. For EV payments,  percentages reflect a cost reduction for EV owners, which is a favorable outcome. A dash (—) indicates cases where the comparison is undefined due to the baseline value being zero (e.g., no export occurred in the flat-rate or charge-only approaches).

\begin{table}[H]
\centering
\caption{Improvement of proposed V2G menu relative to each benchmark scheme.}\vspace{-2mm}
\label{tab:summary_improvement}
\begin{tabularx}{\linewidth}{l
                 >{\centering\arraybackslash}X
                 >{\centering\arraybackslash}X
                 >{\centering\arraybackslash}X}
\toprule
\textbf{Compared to} &
\textbf{Increase in Exported energy (\%)} &
\textbf{Increase in Parking-lot profit (\%)} &
\textbf{Reduction in EV payments (\%)} \\
\midrule
Baseline 1 (Adj.\ RT)   & 87.3\% & 29.61\% & 13.7\% \\
Baseline 2 (Flat)       &  —       & 22.91\% & 18.47\% \\
Baseline 3 (Hybrid)     & 235.72\% & 25.97\% & 9.59\% \\
Charge-only Approach    &  —       & 29.61\% & 17.01\% \\
\bottomrule
\end{tabularx}
\end{table}
\section{Discussion}
\label{sec:discussion}

The proposed dynamic menu-based pricing mechanism is designed for real-time operation and can be implemented using data already available in commercial charging settings. Its rolling-horizon structure enables smooth updates as new vehicles arrive, making the approach practical for deployment.

A key observation from the study is the robustness of the framework across a wide range of conditions. The mechanism maintains robust behavior under variations in feeder capacity, menu granularity and seasonal wholesale price profiles. The sensitivity analysis further shows that substantial price forecast noise does not significantly shift operator outcomes, indicating that the mechanism is robust even in environments with imperfect or volatile market information.

The modular structure of the formulation also supports broader extensions. Heterogeneous user preferences, additional distributed energy resources or coordination across multiple parking sites can be incorporated within the optimization architecture. These extensions would expand applicability without altering the core mechanism.
\section{Conclusion}
\label{sec:conclusion}
We presented a dynamic menu-based pricing mechanism for EV charging in a parking lot setting that also enables EVs to discharge energy back to the grid. The mechanism sets charging and discharging prices for each EV based on its declared needs and real-time wholesale market prices. By aligning operator incentives with user flexibility, it helps parking lot operators increase profit, lowers costs for EV owners, and enhances grid support compared with widely used tariff baselines. Operating in a rolling-horizon mode, the approach updates all schedules whenever a new EV arrives, ensuring real-time responsiveness without reliance on forecasts.

Our results demonstrate that the proposed menu-based pricing is both practical and effective for managing EV charging and V2G. The mechanism remains robust under variations in feeder capacity, fleet size, seasonal price patterns, menu granularity and price fluctuation, highlighting its scalability and adaptability.

\bibliographystyle{IEEEtran}
\bibliography{bib}
\end{document}